\documentclass[a4paper,12pt, reqno]{amsart}
\usepackage{mathrsfs,amsmath,amssymb,latexsym,amsfonts, amscd}
\usepackage[all,ps,cmtip]{xy}%\setlength{\parindent}{.4 in}
\let\mathcal\mathscr

\def\llra{\hbox to 10mm{\rightarrowfill}}

\def\lllra{\hbox to 15mm{\rightarrowfill}}

\def\phi{{\varphi}}

\def\cO{\mathcal{O}}

\DeclareMathOperator{\Pic}{Pic}

\DeclareMathOperator{\id}{id}

\newcommand{\bC}{{\mathbb C}}

\newcommand\lrw{\longrightarrow}

\newcommand\OO{{\mathcal{O}}}

\newtheorem{lemm}{Lemma}[section]
\newtheorem{theo}[lemm]{Theorem}

\newtheorem*{conj*}{Conjecture}

\theoremstyle{definition}

\newtheorem{rema}[lemm]{Remark}

\newtheorem{exam}[lemm]{Example}

\newtheorem{qu}[lemm]{Question}

\theoremstyle{remark}
\newtheorem*{remark*}{Remark}
\newtheorem*{note*}{Note}

\title[Unboundedness of fiber invariants]{Unboundedness of fiber invariants of canonically fibred varieties of general type}
\author{Meng Chen}
\address{\rm School of Mathematical Sciences \& LMNS, Fudan University, Shanghai 200433, China}\email{mchen@fudan.edu.cn}

\author{Zhi Jiang}
\address{\rm D\'epartement de Math\'ematiques, B\^{a}timent 425,
Universit\'{e} Paris-Sud,
F-91405 Orsay, France} \email{zhi.jiang@math.u-psud.fr}

\thanks{The first author was supported by National Natural Science Foundation of China (\#11231003, \#11171068, \#11121101)}
\begin{document}
\begin{abstract} We answer an open question concerning the boundedness of canonical fiber spaces in high dimensions and prove the following:  for any set of integers $n\geq 3$, $0<d<n$ and $N>0$,  there exists a nonsingular projective $n$-fold $X$ of general type so that $X$ is canonically fibred by $d$-dimensional varieties $F$ with $p_g(F)\geq N$.  This  disproves the desired boundedness parallel to Beauville's boundedness theorem (\cite{Bea}) in the surface case.
\end{abstract}
 
\maketitle

\section{\bf Introduction}
We work over any algebraically closed field $k$ of characterictic $0$ (for instance, $k=\bC$). In this article, we assume that $X$ is a minimal projective variety of general type,  namely the canonical (Weil) divisor $K_X$ is nef and big and $X$ has at worst ${\mathbb Q}$-factorial terminal singularities. We call the integer $p_g(X):=\dim_k H^0(X, K_X)$ {\it the geometric genus of $X$}. Denote by $|K_X|$ the canonical linear system of $X$.  
The smallest positive integer $r_X$ such that $r_XK_X$ is Cartier is called {\it the Cartier index of $X$}. When $r_X=1$, we say that $X$ is {\it Gorenstein}. 

 Assume $p_g(X)\geq 2$, the global sections of ${\mathcal O}_X(K_X)$ defines a map $$\varphi_1:=\Phi_{|K_X|}:X\dashrightarrow \mathbb{P}^{p_g(X)-1},$$ which is called {\it the canonical map of $X$}.  
The rational map $\varphi_1$ belongs to two types: (I) $\varphi_1$ is generically finite; (II) $\dim\overline{\varphi_1(X)}<\dim X$, in which case we say that $X$ {\it is canonically fibred by $d$-dimensional varieties}, where $d:=\dim X-\dim\overline{\varphi_1(X)}$.

%Assume that $P_g(X)\geq 2$. Then we have the canonical map
%$$\phi_{K_X}: X\dashrightarrow \mathbb{P}^{P_g(X)-1}.$$
%Denote by $W$ the image of $\phi_{K_X}$. We say that $X$ is canonically of fiber type if $0<\dim W<\dim X$.

When $n=1$, it is well-known that $\varphi_1$ is either birational or a double covering.

When $n=2$, it is clear from the surface theory (see, for instance, \cite{Bea, BHPV}) that either $\varphi_1$ is generically finite of degree $\leq 36$ or $X$ is canonically fibred by curves of genus $\leq 36$.

The natural and interesting question in high dimensions is whether there is the  ``canonical boundedness'' parallel to the cases $n\leq 2$. Unfortunately, when $\varphi_1$ is generically finite, the generic degree $\deg(\varphi_1)$ can be arbitrarily large due to the interesting example of Hacon \cite[Example 2.1]{hac} in dimension 3. Nevertheless, to the authors's knowledge, the following question remained open:

\begin{qu}\label{question} Let $n\geq 3$. Does there exist a constant $N(n)>0$ so that, for all canonically fibred smooth projective $n$-folds $X$ of general type, any smooth model $F$ of an irreducible component of a general fiber of $\varphi_1$ has the geometric genus $p_g(F)<N(n)$?
%For canonically fibred varieties of general type of dimension $n$, are the fibers $F$ birationally bounded? More precisely, we may ask if there exists a positive number $N(n)$ depending only $n$, such that $P_g(F)\leq N(n)$ holds true for any canonically fibred variety $X$?
\end{qu}

The aim of this note is then to give the negative answer to Question \ref{question}.  To state the main result, we first fix necessary notations.  Assume $X$ is canonically of fiber type.  Take a birational modification $\pi:X'\lrw X$ so that $X'$ is smooth and $\varphi_{|K_{X'}|}$ is a morphism.   {}From the Stein factorization of
$\varphi_{|K_{X'}|}$, one gets an induced fibration $f:X'\lrw \Gamma$ onto a lower dimensional projective variety $\Gamma$.  Denote by $F$ a general fiber of $f$.  Our main result is the following:

\begin{theo}\label{main}  For any set of integers $n\geq 3$, $0<d<n$ and $N>0$, there exists a nonsingular projective $n$-dimensional variety $X$ of general type so that $X$ is canonically fibred by $d$-dimensional sub-varieties $F$ with $p_g(F)>N$.
\end{theo}

Theorem \ref{main}, implied by Examples \ref{E1}, \ref{E2} and Example \ref{E3},  clearly disproves the desired ``canonical boundedness'' in high dimensions.  In fact, Example \ref{E1} does the case $d=n-2$. Example \ref{E2} does the cases: $1\leq d\leq n-2$ and Example \ref{E3} does the case $d=n-1$.  
\bigskip

 %Taking the Stein factorization of $g$, we get $g: X'\xrightarrow{f}\Gamma\xrightarrow{s} W$, where $\Gamma$ is normal projective and $f$ is a fibration. Denote by $F$ a general fiber of $f$. The following question has been asked:
%Beauville proved that for a canonically fibred surface $X$, if $\chi(\omega_X)\geq 21$, then $2\leq g(F)\leq 5$ (see \cite[Proposition 2.1]{Bea}).
%This gives an affirmative answer to Question \ref{question} in dimension $2$.

%In this note, we construct several examples to show that there exist minimal varieties (with terminal singularities) of general type, canonically of fiber type with the geometric genus of a general fiber as large as possible. Hence, the answer to Question \ref{question} is negative.

Throughout this paper, the symbol ``$X_{n,d}$'' means a projective minimal (terminal) $n$-fold which is canonically fibred by sub-varieties of dimension $d$ ($1\leq d<n$).

\section{\bf Construction of examples}

%Example \ref{E1} does the case $d = n-2$. Here every smooth fibre $Y$ is a fixed double cover of an Abelian $d$-fold.

\begin{exam}\label{E1} {\bf The variety $X_{n,n-2}$.} We consider a nonsingular projective minimal surface $S$, with $p_g(S)=4$ and $q(S)=0$, which was first constructed by Beauville \cite[Proposition 3.6]{Bea}. We recall that, such a surface $S$ has an involution $\sigma_S$ of which the fixed locus consists of $20$ isolated points so that the quotient surface $T:=S/\langle\sigma_S\rangle$ is a quintic surface (with only $20$ ordinary double points) in $\mathbb{P}^3$. Moreover, by the construction, the canonical system $|K_S|$ is base-point-free and is hence $\sigma_S$-invariant. Denote by $\pi: S\rightarrow T$ the corresponding double cover. Then one has $|K_S|=\pi^*|K_T|$ or, in other words, $\Phi_{|K_S|}$ factors through $\pi$ and $\Phi_{|K_T|}$.
%$H^0(S, \omega_S)=\pi^*H^0(T, \omega_T)$.

Now we take any double covering $Y\rightarrow A$ from a smooth projective variety of general type to an abelian variety of dimension $d>0$ with the corresponding involution $\sigma_Y$. Define $X$ to be the quotient $$X_{d+2,d}:=(Y\times S)/\langle\sigma_Y\times\sigma_S\rangle.$$ Then $X$ has  terminal singularities and, by considering the natural quotient $f: X\rightarrow A\times T$, we see that
$|K_X|=f^*|K_{A\times T}|$ and that $p_g(X)=h^0(A\times T,\omega_{A\times T})$. Hence
%$\omega_X$ is globally generated and
$\Phi_{|K_X|}$ naturally factors through the following quotient morphism:
$$X:=(Y\times S)/\langle\sigma_Y\times\sigma_S\rangle\rightarrow S/\langle\sigma_S\rangle=T$$
and $\Phi_{|K_T|}$. Since $n:=\dim X=d+2$ and $\dim\Phi_{|K_T|}(T)=2$, we see that $X$ is canonically fibred by varieties (isomorphic to $Y$) of dimension $d$. It is clear that $p_g(Y)$ can be arbitrary large.  
%We see that $Y$ is a general fiber of the canonical map of $X$. We can choose $Y$ with $P_g(Y)$ as large as possible and in this way, we see the unboundedness of fibers of canonical maps.
\end{exam}

{\em {}From the construction of $X_{n,n-2}$, we see that the induced fibration from $\Phi_{|K_{X_{n, n-2}}|}$ is iso-trivial with the general fiber a fixed double cover of the given  abelian $d$-fold $A$.}

%In the above example, we get $n$-folds canonically fibred by constant fibers of dimension $n-2$ with arbitrary large geomeric genus.
%We can mimic a construction of Hacon \cite{hac} to get examples canonically fibred by varieties with nonconstant moduli and arbitrary large geomeric genus.

\begin{exam} \label{E2}{\bf Varieties $X_{n,d}$ with $1\leq d\leq n-2$}.  Let $A_i$ be a simple principally polarized abelian variety of dimension $d_i>0$ with a principal polarization $\Theta_i$ for $i=1,2$.  The existence of such pairs $(A_i, \Theta_i)$ is well-known. 
Let $\rho_i: F_i\lrw A_i$ be the double cover ramified over a general member of $|2\Theta_i|$. Hence $F_i$ is a smooth variety of general type.  Write
\begin{equation}\label{e1}\rho_{i*}\omega_{F_i}=\cO_{F_i}\oplus\Theta_i
\end{equation}
and denote by $\sigma_i$ the involution corresponding to $\rho_i$ for each $i=1,2$. Clearly we have $p_g(F_i)=2$ for each $i$. 

Let $S$ be a smooth projective variety of dimension $d_3\geq 1$ with $p_g(S)=0$. Take a very ample line bundle $L$ on $S$  such that
$K_S + L$ is also very ample. Pick a smooth divisor $D\in |2L|$.
Let $\rho_3: T\rightarrow S$ be a double cover of $S$ ramified in $D$. Then $T$ is a smooth projective variety of general type and we have
\begin{eqnarray}\label{e2}\rho_{3*}\omega_{T}&=&\omega_{S}\oplus (\omega_S\otimes L). 
\end{eqnarray} Denote by $\sigma_3$ the   involution  corresponding to $\rho_3$.

{\it Step 1}. We consider the product $F_1\times F_2\times T$ and the involution $\sigma=\sigma_1\times\sigma_2\times \sigma_3$.  Define $X$ to be the quotient $(F_1\times F_2\times T)/\langle \sigma\rangle$. Then $X$ is a minimal projective variety with terminal singularities of Cartier index 2 and $K_X^3>0$.

Let $\rho: X\rightarrow Z:=A_1\times A_2\times S$ be the natural morphism which is clearly a $(\mathbf{Z}/2\mathbf{Z}\times \mathbf{Z}/2\mathbf{Z})$-covering. Moreover, by equations (\ref{e1}) and (\ref{e2}), we have
 \begin{eqnarray}\label{pushforward}
\nonumber\rho_*\omega_X &=& \big(\cO_{A_1}\boxtimes \cO_{A_2}\boxtimes\omega_S\big)\oplus \big(\Theta_1\boxtimes \Theta_2\boxtimes \omega_S\big)\\ && \oplus  \big(\Theta_1\boxtimes \cO_{A_2}\boxtimes(\omega_S\otimes L)\big)\oplus \big(\cO_{A_1}\boxtimes \Theta_2\boxtimes(\omega_S\otimes L)\big),
\end{eqnarray}
where $\cO_{A_1}\boxtimes \cO_{A_2}\boxtimes \omega_S:=\iota_1^*\cO_{A_1}\otimes \iota_2^*\cO_{A_2}
\otimes \iota_3^*\omega_S$ and $\iota_i$
is the projection from $A_1\times A_2\times S$
onto the $i$-th factor for $i=1,2,3$.

In particular, we have 
\begin{eqnarray*}H^0(X, \omega_X)&\cong & H^0(Z, \cO_{A_1}\boxtimes \cO_{A_2}\boxtimes\omega_S)\oplus H^0(Z, \Theta_1\boxtimes \Theta_2\boxtimes \omega_S)\\
&&\oplus H^0(Z, \Theta_1\boxtimes \cO_{A_2}\boxtimes (\omega_S\otimes L))
\oplus H^0(Z, \cO_{A_1}\boxtimes \Theta_2\boxtimes (\omega_S\otimes L)).
\end{eqnarray*}
Besides, we have the following commutative diagram:
\begin{eqnarray*}
\xymatrix{
X\ar[dr]_f\ar[r]^{\rho} & Z\ar[d]^{\iota_3}\\
&T/\langle\sigma_3\rangle=S,}
\end{eqnarray*}
where $f:=\iota_3\circ \rho$. By (\ref{pushforward}), we have 
\begin{eqnarray*}f_*\omega_X\simeq &&\omega_S\oplus \omega_S \oplus(\omega_S\otimes L)\oplus  (\omega_S\otimes L).
\end{eqnarray*}
For the general point $t\in S$, denote by $X_t$ the general fiber of $f$ over $t$.  By direct calculation, one has $X_t\simeq F_1\times F_2$ and, clearly, $p_g(X_t)=p_g(F_1)\cdot p_g(F_2)=4$. Then the restriction map $$H^0(X, \omega_X)\rightarrow H^0(X_t, \omega_{X_t})$$ has the image dimension $=2$ since $h^0(A_i, \Theta_i)=1$ for $i=1,2$ by assumption. 
Here we note that $f_*\omega_X$ has rank $4$.
Hence, the image $\overline{\Phi_{|K_X|}(X)}$ is birational to a $\mathbf{P}^1$-bundle $V$ over $S$. Thus we have the following diagram:
\begin{eqnarray*}
\xymatrix{
X\ar@{.>}[r]^{\Phi_{|K_X|}}\ar[dr]_f& V\ar[d]\\
&S.}
\end{eqnarray*}

Denote by $F$ an irreducible component of a general fiber of $\Phi_{|K_X|}$. As we have already seen, $F$ is a sub-canonical divisor of $K_{X_t}=K_{F_1\times F_2}$. In fact, $F$ moves in a sub-pencil of $|K_{X_t}|$.
%Then $F$ is  a canonical divisor of a general fiber $F_1\times F_2$ of $f$.  It is easy to see that $F$ is irreducible.
\medskip

{\it Step 2}. We note that the natural morphism $\alpha: X\rightarrow A_1\times A_2$ is the Albanese morphism of $X$. We take $P=P_1\boxtimes P_2$, a general $m$-torsion line bundle of $X$ and let $\pi_m: X_m\lrw X$ be the \'etale $m$-cover of $X$ induced by $P$. Then we have $\pi_{m*}\omega_{X_m}=\oplus_{i=0}^m (\omega_X\otimes P^{\otimes i})$. Since both $P_1$ and $P_2$ are of order $m$, by (\ref{pushforward}), we have $H^0(X_m, \omega_{X_m})\cong H^0(X, \omega_X)$. Thus, we have the following commutative diagram:
\begin{eqnarray*}
\xymatrix{
X_m\ar[r]^{\pi_m} \ar@{.>}[dr]_{\Phi_{K_{|X_m|}}}& X\ar@{.>}[d]^{\Phi_{|K_X|}}\\
& V}
\end{eqnarray*}
Moreover, since $\alpha(F)$ generates $A_1\times A_2$, we can choose a very general $m$-torsion $P$ such that $P|_F$ is also a $m$-torsion (noting that $F$ moves in an algebraic family on $X_t$). Thus, over $F$, the irreducible component $F_m$ of the general fiber of $\Phi_{|K_{X_m}|}$ is simply an  \'etale $m$-cover of $F$.

Take a desingularization $\varepsilon: F'\rightarrow F$.  Then  the induced \'etale cover $F_m'$ of $F'$ is a desingularization of $F_m$.  Consider the induced generically finite morphism $\beta:=\alpha|_F\circ \varepsilon: F'\rightarrow A_1\times A_2$.
Since $F'$ is of general type and both $A_1$ and $A_2$ are simple, we have
$H^0(F', \omega_{F'}\otimes \beta^*Q)\neq 0$, for any $Q\in \Pic^0(A_1\times A_2)$ (see \cite[Corollary 3.5]{CDJ}).
 Therefore, $\lim_{m\rightarrow \infty }p_g(F_m')\rightarrow+\infty$.

We see that $X_m$  is minimal projective with terminal singularities and that  $X_m$ is canonically fibred by  varieties which are algebraically equivalent to $F_m$. Moreover,  $\dim X_m=d_1+d_2+d_3\geq 3$ and $\dim F_m=d_1+d_2-1$.   Hence, choosing appropriate $d_1$, $d_2$, and $d_3$,  $\dim F_m$ can be any integer between $1$ and $\dim X_m-2$. We set $X_{n,d}^{<m>}:=X_m$ where $n=d_1+d_2+d_3$ and $d=d_1+d_2-1$.
\end{exam}

{\em Clearly, as $m$ goes to infinity, our variety  $X_{n,d}^{<m>}$ is a minimal $n$-fold of general type which is canonically fibred by sub-varieties $F_m$ of dimension $d$ with $1\leq d\leq n-2$ and, furthermore,  the smooth model of $F_m$ has arbitrarily large geometric genus. }

%In  Example \ref{E1}, the canonical map is isotrivial. We now construct varieties canonically fibred by varieties with nonconstant moduli, where $1\leq d\leq n-2$ and a general fiber $F$ is a $(\mathbf{Z}/2\mathbf{Z}\times\mathbf{Z}/2\mathbf{Z})$-cover of abelian $d$-fold.
%Example \ref{E3} does the case $d=n-1$.

\begin{exam}\label{E3}{\bf The variety $X_{n, n-1}$}.  We consider a variant of Chen--Hacon's example \cite[\S 4]{jch} where a threefold  of general type, with maximal Albanese dimension and with $p_g=1$, was constructed. 

{\it Step 1}. To be precise, let $(A_i, \Theta_i)$ be simple principal polarized abelian varieites for $i=1,2,3$. Let $\rho_i: F_i\rightarrow A_i$ be the double covering ramified over a general element of $|2\Theta_i|$, with the induced involution $\sigma_i$, for each $i$. Choose $P_i\in \Pic^0(A_i)$ of order $2$, we consider the induced \'etale covers $F_i'\rightarrow F_i$, with the induced involution $\iota_i$. The involution $\sigma_i$ on $F_i$ pulls back to an involution $\sigma_i'$ on $F_i'$.

 We then define $X$ to be the quotient of $F_1'\times F_2'\times F_3'$ by the group of automorphisms generated by $\id_{F_1'}\times \sigma_2'\times \iota_3$, $\iota_1\times \id_{F_2'}\times \sigma_3'$, $\sigma_1'\times\iota_2\times \id_{F_3'}$, $\sigma_1'\times \sigma_2'\times \sigma_3'$.

One can check that $X$ is minimal projective of general type with terminal singularities of Cartier index $2$. Moreover, we have a $(\mathbf{Z}/2\mathbf{Z}\times \mathbf{Z}/2\mathbf{Z})$-cover:
$$\rho: X\rightarrow A:=A_1\times A_2\times A_3,$$
and by analyzing  the  action of the Galois group, we have
\begin{eqnarray}\label{pushforward2}\nonumber \rho_*\omega_X&=&\cO_{A}\oplus \big(\Theta_1\boxtimes (\Theta_2\otimes P_2)\boxtimes P_3\big)\oplus\big((\Theta_1\otimes P_1)\boxtimes P_2\boxtimes \Theta_3\big)\\&&\oplus \big(P_1\boxtimes \Theta_2\boxtimes (\Theta_3\otimes P_3)\big).
\end{eqnarray}
We note that $\rho$ is the Albanese morphism of $X$.
Modulo one more \'etale cover of $A_1$ to annihilate $P_1$ upstairs,
%After taking \'etale covers of $A_1$,
we may and do assume that $P_1=\cO_{A_1}$ is trivial. In other words, we have got such a variety $X$ that $\rho_*\omega_X$ has exactly two summands contributing to its global sections, which are $\cO_{A}$ and $\OO_{A_1}\boxtimes \Theta_2\boxtimes (\Theta_3\otimes P_3)$ by (\ref{pushforward2}).  Hence
we have $\dim H^0(X, \omega_X)=2$. We denote $\dim X=n\geq 3$,  then automatically $X$ is canonically fibred by varieties of dimension $n-1$.

{\it Step 2}. We consider the the natural morphism $f: X\rightarrow A_1$. Then $f_*\omega_X=\cO_{A_1}\oplus\cO_{A_1}$. Hence, the restriction of the canonical pencil $|K_X|$ to a general fiber of $f$ is still a pencil. Then, up to birational equivalence, $\Phi_{|K_X|}$ does not factor through $f$. Let $F$ be any irreducible component of the general fiber of $\Phi_{|K_X|}$. Then $f|_F: F\rightarrow A_1$ is surjective. In particular, for any integer $m>1$, there exists a torsion line bundle $P_m$ on $A_1$ such that $f^*P_m|_{F}$ is still a $m$-torsion.

We then take a desingularization $\varepsilon: F'\rightarrow F$. Considering the induced generically finite morphism: $\beta:=\rho|_F\circ \varepsilon: F'\rightarrow A$,  since $A$ has only $3$ simple factors and $\beta$ is not surjective, we have $H^0(F', \omega_{F'}\otimes \beta^*Q)\neq 0$, for any $Q\in \Pic^0(A)$ (see \cite[Proposition 4.3]{CDJ}).

{}Finally let $\pi_m: X_m\rightarrow X$ be the \'etale $m$-cover induced by $f^*P_m$. By (\ref{pushforward2}), we can easily check that $H^0(X_m, \omega_{X_m})\cong H^0(X, \omega_X)$ and $\omega_{X_m}=\pi_m^*\omega_X$. Thus, we have the commutative diagram:
\begin{eqnarray*}
\xymatrix{
X_m\ar[r]^{\pi_m} \ar@{.>}[dr]_{\Phi_{|K_{X_m}|}}& X\ar@{.>}[d]^{\Phi_{|K_X|}}\\
& \mathbf{P}^1.}
\end{eqnarray*}
Clearly $|K_{X_m}|$ is composed of a pencil of $(n-1)$-dimensional varieties. Any irreducible element $F_m$ of the general fiber of $\Phi_{|K_{X_m}|}$ is then an \'etale $m$-cover of $F$ induced by $f^*P_m|_{F}$.
Let $F_m'$ be the corresponding \'etale cover of $F'$, then $F_m'$ is a desingularization of $F'$ and simple calculations tell us that $\lim_{m\rightarrow \infty }p_g(F_m')\rightarrow +\infty$. Set $X_{n,n-1}^{<m>}$, which is what we want to realize. 
\end{exam}

{\it Clearly, as $m$ goes to infinity, $X_{n,n-1}^{<m>}$ is a projective minimal (terminal) $n$-fold and is canonically fibred by $(n-1)$-folds $F_m$. Furthermore, a smooth model of $F_m$ has arbitrarily large geometric genus.}

\begin{rema} Though the ``canonical boundedness'' does not hold in  dimension $\geq 3$ by Theorem \ref{main}, we still have boundedness results in special cases. For example, if we consider canonically fibred Gorenstein minimal 3-folds of general type, Chen and Hacon \cite[Theorem 1.1]{CH} showed that the fiber invariants of $\Phi_{|K_X|}$ are upper bounded, and some interesting examples with large fibre invariants were found by Chen-Cui in \cite{CC}.
\end{rema}

It is then natural to ask the following questions:

\begin{qu} Let $m>1$ and $n\geq 3$ be two integers. Considering all those projective minimal $n$-folds of general type which are $m$-canonically fibred, are the birational invariants of fibers of $\Phi_{|mK|}$ universally bounded from above?
\end{qu}

\begin{qu} Let $Y$ be a projective Gorenstein minimal $n$-fold ($n\geq 4$) of general type with $p_g(Y)\geq 2$.
\begin{itemize}
\item[(1)] Assume that $\varphi_1$ is generically finite. Is $\deg(\varphi_1)$ bounded from above?

\item[(2)] Assume that $\varphi_1$ is of fiber type. Is the geometric genus of any smooth model of any irreducible component in the general fiber of $\varphi_1$ bounded from above?
\item[(3)] Is there an effective inequality of the form $K_Y^n\leq c(n)p_g(Y)$ where $c(n)$ is a constant depending only on $n$?
\end{itemize}
%We note that all examples we are considering have essentially singularities of type $\frac{1}{2}(1,1,1)$  and hence are not Gorenstein. We then ask:  for canonically fibred Gorenstein (or even smooth) minimal $n$-fold, are the canonical fibers bounded ?
\end{qu}

\end{document}